\documentclass[11pt]{amsart}
\usepackage{amssymb}
\usepackage{amscd}
\usepackage{amsmath}
\usepackage{amsfonts}
\usepackage{pstricks}
\usepackage{pst-node}
\usepackage{graphicx}
\theoremstyle{plain}
\newtheorem{thm}{Theorem}[section]

\theoremstyle{remark}

\def\pmc#1{\setbox0=\hbox{#1}
    \kern-.1em\copy0\kern-\wd0
    \kern.1em\copy0\kern-\wd0}

\begin{document}

\bigskip

\title[On generalization of Homotopy Axiom]{
On generalization of Homotopy Axiom}

\bigskip

\author[U.~H.~Karimov]{Umed H. Karimov}
\address{Institute of Mathematics,
Academy of Sciences of Tajikistan, Ul. Ainy $299^A$,
Dushanbe 734063,
Tajikistan}
\email{u.h.karimov@gmail.com}


\begin{abstract}
In the paper \cite{Si} it was proven that if $G$ is compact topological group or field then in the Homotopy Axiom for Alexander-Spanier-Kolmogoroff cohomology the parameter segment $[0, 1]$ can be replaced by any compact connected space $T$. The purpose of the paper is to shows that the parameter space $T$ can not be replaced in general by locally compact connected space.
\end{abstract}

\date{\today}

\maketitle

By $H^*$ we mean Alexander-Spanier-Kolmogoroff cohomology or, equivalently for paracompact spaces, Aleksandroff-\v{C}ech cohomology,

The Homotopy Axiom of Eilenberg and Steenrod \cite{ES} for cohomology is equivalent to the following statement:

\medskip

\noindent {\bf Axiom.}   {\it Let $i_0$ and $i_1$ are embeddings of the pair $(X,A)$ to $(X,A)\times [0, 1]$ defined by $i_0(x) = (x,0)$ and $i_1(x) = (x,1)$ respectively, then $i_0^{\ast} = i_1^{\ast}: H^{\ast}(X\times [0, 1] ; G)\rightarrow H^{\ast}(X; G) $.}

\medskip

In the paper \cite{S} was proven the theorem:

\medskip

\begin{thm} Suppose that $X$ and $Y$ are any spaces, $T$ is a compact,
connected space, and $G$ is a compact abelian topological group or a
finite-dimensional vector space. Then for any continuous function
$F: X \times T \rightarrow Y$, one has
$$F_r^{\ast} = F_s^{\ast}:H^{\ast}(Y; G)\rightarrow H^{\ast}(X; G)$$
for all $r, s$ in $T$, where $F_t:X\rightarrow Y$ is defined by $F_t(x) = F(x, t)$.
\end{thm}

The purpose of the paper is to prove following theorem:
\medskip

\begin{thm}
 There exist locally compact connected spaces $X$ and $T$, two points $r$ and $s$ in $T$ such that the embeddings $F_r: X \rightarrow X\times T$ and $F_s: X \rightarrow T\times T$ defined by $F_r (x) = (x, r)$ and $F_s (x) = (x, s)$ respectively induce different homomorphisms in cohomology groups.
\end{thm}

Space $X$ can not be compact space \cite{S}, \cite{W}, the parameter space $T$ can not be acyclic space \cite{Br}.
Let $G = Z_p$ be the field of the characteristic $p, p$ is prime number and
$$T = \{(x,y): (x,y) \in {R^2}, y = \sin (1/x)\ \ \ \text {if} \ \ \ 0 < x \leq 1, \ \ \ \text{and}\ \ \ y\in [-1,0)\ \ \ \text {if} \ \ \ x = 0\}.$$

Let $r = (0, -1)$ and $s = (1, \sin 1).$ Let the space $X$ be the same space as $T$. Consider the mappings $F^1_t: H^1(T \times T; Z_p) \rightarrow H^1(T; Z_p)$ induced by the embeddings $F_t: T \rightarrow T\times \{t\} \subset T\times T.$ Let $V = T \cap \{(x,y) \in R^2, y < 0\}$ and $V_i = V\cup (T\cap \{(x,y) \in R^2, x > \frac{1}{i}\})$. It is obvious that $T = \cup V_i, V_i \subset V_{i+1}$ and $V_i$ are open in $T.$ We have commutative diagram \cite[Theorem 1]{Kh}:
\begin{align}
0 &\rightarrow &  &\underleftarrow{\lim}^{(1)} H^0(V_i\times V_i; Z_p)& &\rightarrow & &H^1(T \times T; Z_p)&  & \rightarrow &  &\underleftarrow{\lim}\ H^1(V_i\times V_i;Z_p)& \rightarrow 0\nonumber\\
& & & f_r \downarrow \downarrow f_s & &  & & F_r^1\downarrow \downarrow F_s^1 & & & & \downarrow & \nonumber \\
0 &\rightarrow & &\underleftarrow{\lim}^{(1)}\ H^0(V_i; Z_p)&  &\rightarrow &  &H^1(V_i; Z_p)&  & \rightarrow & &\underleftarrow{\lim} H^1(V_i; Z_p)& \rightarrow 0 \nonumber
\end{align}

Since $\underleftarrow{\lim}H^1(V_i\times V_i; Z_p) = 0$ and $\underleftarrow{\lim}H^1(V_i; Z_p) = 0$ it is possible to identify corresponding groups and identify $f_r$ with $F_r^1$ and identify $f_s$ with $F_s^1$ respectively where the mappings $f_r$ and $f_s$ are generated by the embeddings $ V_i \subset V_i\times r \subset V_i\times V_i$ and $ V_i \subset V_i\times s \subset V_i\times V_i$ respectively. The group $\underleftarrow{\lim}^{(1)}H^0(V_i; Z_p)$ is nontrivial because the groups $H^0(V_i; Z_p)$ are countable and the inverse spectrum $\{H^0(V_i; Z_p)\}$ does not satisfy the Mittag-Leffler condition, see e.g. \cite[Proposition 1]{Kh}.

Consider commutative diagram:
\begin{align}
 &\underleftarrow{\lim}^{(1)} H^0(V_i\times V_i; Z_p)& &\stackrel{h}{\longrightarrow} & &\underleftarrow{\lim}^{(1)} H^0(V_i\times \{r,s\}; Z_p)&  \nonumber\\
 & f_r \downarrow \downarrow f_g & &  & &\downarrow \textit{id} &  \nonumber \\
 &\underleftarrow{\lim}^{(1)}\ H^0(V_i; Z_p)&  &\stackrel{\phi,\psi}{\underleftarrow{\longleftarrow}} &  &\underleftarrow{\lim}^{(1)} H^0(V_i\times \{r,s\}; Z_p)&  \nonumber
\end{align}

The mapping $h$ is generated by the natural embeddings $V_i\times \{r,s\} \subset V_i\times V_i$, \textit{id} is identity mapping and mappings $\phi, \psi$ are generated by the embeddings $V_i$ in $V_i\times \{r,s\}$. From the exact sequence of inverse spectra
$$0 \rightarrow H^0(V_i\times V_i, V_i\times \{r, s\}; Z_p) \rightarrow H^0(V_i\times V_i; Z_p) \rightarrow H^0(V_i\times \{r, s\}; Z_p) \rightarrow 0$$
\noindent it follows that $h$ is epimorphism because for any countable inverse spectrum $\underleftarrow{\lim^{(2)}} = 0.$ Since $\phi\neq \psi$ it follows that $f_r \neq f_s$ and therefore $F_r^1 \neq F_s^1.$ The Theorem is proved.

\end{document}